\documentclass[journal,onecolumn]{IEEEtran}

\usepackage{graphicx}
\usepackage{float}
\usepackage{amsmath}
\usepackage{amssymb}
\usepackage{subcaption}
\usepackage{array}
\usepackage{colortbl}
\usepackage[T1]{fontenc}
\usepackage{xstring}
\usepackage{setspace}
\usepackage{xspace}
\usepackage{xcolor}
\usepackage{algorithm}
\usepackage{algpseudocode}
\usepackage{longtable}
\usepackage{booktabs}
\usepackage{url}
\usepackage[noadjust]{cite}



\algnewcommand\algorithmicinput{\textbf{Input:}}
\algnewcommand\Input{\item[\algorithmicinput]}
\algnewcommand\algorithmicinputcont{\phantom{Input: }}
\algnewcommand\Inputcont{\item[\algorithmicinputcont]}

\newcommand{\Ilambda}{I_{\lambda}}

\newcommand{\lmin}{\underline{\lambda}} 
\newcommand{\lmax}{\overline{\lambda}}

\newcommand{\gmin}{\lambda_{\textrm{min}}} 
\newcommand{\gmax}{\lambda_{\textrm{max}}} 

\newcommand{\Yblock}{Y} 
\newcommand{\U}{U} 
\newcommand{\Xapprox}{\tilde{X}} 
\newcommand{\Xexact}{X_{\lambda}} 
\newcommand{\xexact}{\mathbf{x}} 
\newcommand{\xapprox}{\tilde{\mathbf{x}}} 
\newcommand{\Nmat}{n} 
\newcommand{\ewexact}{\lambda} 
\newcommand{\ewapprox}{\tilde{\lambda}} 

\renewcommand{\textsubscript}[1]{\ensuremath{{}_{\mathtt{#1}}}}
\newcommand{\beastmx}{\texttt{BEAST-M\textsubscript{\ast,x}}\xspace} 
\newcommand{\beastm}{\texttt{BEAST-M}\xspace} 
\newcommand{\beastmn}{\texttt{BEAST-M\textsubscript{\ast,n}}\xspace} 
\newcommand{\beastmout}{\texttt{BEAST-M\textsubscript{o,\ast}}\xspace} 
\newcommand{\beastmin}{\texttt{BEAST-M}\textsubscript{i,\ast}\xspace} 
\newcommand{\beastmxout}{\texttt{BEAST-M\textsubscript{o,x}}\xspace} 
\newcommand{\beastmxin}{\texttt{BEAST-M}\textsubscript{i,x}\xspace} 
\newcommand{\beastmnout}{\texttt{BEAST-M\textsubscript{o,n}}\xspace} 
\newcommand{\beastmnin}{\texttt{BEAST-M}\textsubscript{i,n}\xspace} 
\newcommand{\beastc}{\texttt{BEAST-C}\xspace} 
\newcommand{\beast}{\texttt{BEAST}\xspace} 

\newcommand{\beastcad}{\texttt{BEAST-C\textsubscript{ad}}\xspace} 
\newcommand{\beastcn}{\texttt{BEAST-C\textsubscript{n}}\xspace} 
\newcommand{\beastp}{\texttt{BEAST-P}\xspace} 
\newcommand{\nexpect}{\tilde{m}} 
\newcommand{\ntrue}{m} 
\newcommand{\contour}{C} 
\newcommand{\mo}{\mathrm{width}( \U )} 
\newcommand{\quadpts}{q} 
\newcommand{\z}{z} 
\newcommand{\zj}{z_j} 
\newcommand{\wj}{\omega_j} 
\newcommand{\momentsk}{k} 
\newcommand{\nrhs}{\RHSSingle} 
\newcommand{\nmoments}{s} 
\newcommand{\dz}{\mathrm{d}z}
\newcommand{\snc}{\textrm{SNC}}
\newcommand{\sseig}{\texttt{sseig}\xspace}
\newcommand{\strumpack}{\texttt{strumpack}\xspace}
\newcommand{\sseigtobeastc}{\texttt{sseig2beastc}\xspace}
\newcommand{\ghost}{\texttt{ghost}\xspace}

\newcommand{\RHSSingle}{\mathrm{RHS}_{1}}
\newcommand{\RHSovl}{\mathrm{RHS}_{\mathrm{ovl}}}
\newcommand{\BLSovl}{\mathrm{BLS}_{\mathrm{ovl}}}
\newcommand{\scamac}{\texttt{ScaMaC}\xspace}
\renewcommand{\arraystretch}{1.2}

\begin{document}

\title{Flexible subspace iteration with moments for an effective contour integration-based eigensolver}

\author{Sarah~Huber\textsuperscript{1,a}, Yasunori~Futamura\textsuperscript{b}, Martin~Galgon\textsuperscript{a}, Akira~Imakura\textsuperscript{b}, Bruno~Lang\textsuperscript{a}, Tetsuya~Sakurai\textsuperscript{b} \\ \textsuperscript{a} University of Wuppertal, Faculty of Mathematics and Natural Sciences \\ \textsuperscript{b} University of Tsukuba, Department of Computer Science
\thanks{\textsuperscript{1} E-mail: shuber@math.uni-wuppertal.de}
}

\maketitle

\begin{abstract}
Contour integration schemes are a valuable tool for the solution of difficult interior eigenvalue problems. 
However, the solution of many large linear systems with multiple right hand sides may prove a prohibitive computational expense. 
The number of right hand sides, and thus, computational cost may be reduced if the projected subspace is created using multiple moments. 
In this work, we explore heuristics for the choice and application of moments with respect to various other important parameters in a contour integration scheme. 
We provide evidence for the expected performance, accuracy, and robustness of various schemes, showing that good heuristic choices can provide a scheme featuring good properties in all three of these measures.
\end{abstract}

\begin{IEEEkeywords}
Contour-integral based eigensolver; Sakurai--Sugiura methods, subspace iteration
\end{IEEEkeywords}

\section{Introduction}

Contour integration schemes have become popular in recent years as a
broadly applicable spectral filtering approach to interior eigenvalue
problems.
Two of the most well-known of these are Sakurai--Sugiura methods (SSM),
the first of which was introduced by Sakurai and Sugiura in 2003
\cite{sakurai2003projection}, and FEAST, introduced by Polizzi in 2009
\cite{polizzi2009density}.  Both SSM \cite{sakurai2007cirr,relationships2016,
imakura2014block,ikegami2010contour,ikegami2010filter,Imakura2016,doi:10.1260/1748-3018.7.3.249,sakurai2019scalable,imakura2018block} and FEAST \cite{peter2014feast,guettelzolotarev2014,yin2014feast,kestyn2016feast,
gavin2018feast,gavin2018krylov}
have been extended and further developed.
These methods are in particular applicable to the generalized interior
eigenproblem in the positive definite case, i.e., to finding the
eigenpairs $(\ewexact, \xexact)$ of
\begin{equation}
  \label{eqn:genevp}
  A\xexact = B\xexact \ewexact, \qquad
  \ewexact \in \Ilambda = \left[ \lmin, \lmax \right],
\end{equation}
with hermitian $A$ and hermitian positive definite $B$.

The identifying and typically most computationally costly component of
contour integration schemes is the evaluation of a contour integral.  
More precisely, for \eqref{eqn:genevp} one chooses a contour $\contour$
in the complex plane such that $\contour$ contains the eigenvalues in $\Ilambda$, but no
other eigenvalue, and then computes approximations to integrals of the
form
\begin{equation}
  \label{eq:integral}
  \frac{1}{2\pi i}
  \int_{\contour}\z^{\momentsk} (\z B-A)^{-1} B \Yblock\dz
  \; .
\end{equation}
From these, an approximate subspace, $\U$, for the desired eigenspace,
$\Xexact$, is obtained.
In both FEAST and Sakurai--Sugiura Rayleigh--Ritz (SS--RR)
\cite{sakurai2007cirr}, the type of SSM considered in this work, the
approximate eigenvalues and vectors are then obtained from $\U$ with a
Rayleigh--Ritz reduction \cite{saad2003iterative} and solve the reduced eigenproblem
\begin{equation}
  \label{eq:RR}
  A_{\U} W = B_{\U} W \Lambda
\end{equation}

where $ A_{\U} = U^{\ast} A U, B_{\U} = U^{\ast} B U$ and $ \Xapprox := \U W $.

Solvers use numerical quadrature to approximate the integral
\eqref{eq:integral}, acting on a set of initial vectors $\Yblock$ to
construct the approximate subspace $\U$.
Given a quadrature rule with $\quadpts$ quadrature nodes $\zj$ and
coefficients $\omega_j$, $\quadpts$ linear systems of the form
$(\zj B-A)^{-1}B\Yblock$ with multiple right hand sides must be solved.
The defining difference between FEAST and SSM is the use of moments in
the contour integrals.

In Sakurai--Sugiura methods, the approximate integral contains a term
$\z^{\momentsk}$ representing the $\momentsk^{\mathrm{th}}$ moment,
leading to the subspace components
\begin{equation*}
  \U_k = \sum_{j=1}^{\quadpts_{\textrm{SSM}}}
  \wj \zj^{\momentsk} (\zj B-A)^{-1} B\Yblock
  ,
\end{equation*}
which are combined, given a total of $s$ moments, to form the complete approximate subspace as
\begin{equation}
  \label{eqn:momentssubspace}
  \U=[\U_0, \U_1, \U_2, ... ,\U_{\nmoments-1} ].
\end{equation}
In FEAST, only the $0^{\mathrm{th}}$ moment is considered, and $\U$ is
formed as
\begin{equation}
  \label{eqn:contourexact}
  \U \equiv \U_{0}
  = \sum_{j=1}^{\quadpts_{\textrm{FEAST}}}
  \wj (\zj B-A)^{-1} B \Yblock
  .
\end{equation} 

Often the most expensive step of a contour integration scheme is the
solution of linear systems $(\zj B-A)^{-1} B\Yblock$ involved in the
construction of $\U$.
Discounting differences based on the type of linear solver used, the
cost for one integration depends on $q$, the number of quadrature nodes (or block linear systems solved), and on $\RHSSingle$, the number
of right hand sides (or columns of $\Yblock$) in each block linear
system.
The overall cost will therefore be dependant on $\RHSovl$ and
$\BLSovl$, the overall number of right hand sides and of block linear
systems solved, resp., over all potential iterations.

Note that the number of columns of $\U$, $\mo$, is by construction
$\mo=\RHSSingle \cdot \nmoments$, where $\RHSSingle$ is the number of
columns of $\Yblock$, and $\nmoments$ is the number of moments.
Therefore, a subspace $\U$ constructed with $\nmoments$ moments
requires the input matrix $\Yblock$ to have only
$1/\nmoments^{\mathrm{th}}$ of the desired number of columns of $\U$.

SS--RR has traditionally used a large overall subspace size $\mo$,
resulting in rapid convergence rates and thus requiring few iterations.
However, allowing for more iterations can give more possibilities for
adaptive techniques and overall cost reduction.
Predicting the expected behaviour of the algorithms is in this case
necessary to develop good heuristics.
Having previously explored adaptivity with respect to quadrature nodes
\cite{NLA:NLA2124}, we now consider the benefits of allowing the solver
type itself to vary between iterations.
This work considers the use of multiple moments in these schemes, as in
SS--RR, as well as when a single-moment strategy (FEAST) is suitable,
for a given subspace size.
Our goal is to increase the efficiency of the eigensolver by reducing
$\RHSovl$ while maintaining accuracy and robustness.

\section{Flexible subspace iteration}

FEAST is, by definition, a subspace iteration scheme.
Multiple iterations of the sequence
\begin{center}
  \parbox{15cm}{%
    \begin{algorithmic}[1]
      \State Construct subspace $\U$ with contour integration
        according to \eqref{eqn:contourexact}
      \State Rayleigh--Ritz extraction of eigenvalues and vectors
        according to \eqref{eq:RR}
      \State $\Yblock := \Xapprox$
    \end{algorithmic}%
  }
\end{center}
result in a set of vectors, $\Xapprox$, inside (or close to) the desired
eigenspace.

SSM, on the other hand, was not originally defined as an iterative
scheme.
The basic SS--RR roughly works as follows:
\begin{center}
  \parbox{15cm}{%
    \begin{algorithmic}[1]
      \State Construct subspace $\U$ with contour integration and
        moments according to \eqref{eqn:momentssubspace}
      \State Orthogonalize $\U$ and reduce subspace by removing (almost)
        linearly dependent columns
      \State Rayleigh--Ritz extraction of eigenvalues and vectors
        according to \eqref{eq:RR} ,
    \end{algorithmic}%
  }
\end{center}
without iteration.
More recently, the ``inner iteration'' has been introduced \cite{Imakura2016},
\begin{center}
  \parbox{15cm}{%
    \begin{algorithmic}[1]
      \State $\U_{0} := \sum_{j=1}^{\quadpts_{\textrm{SSM}}}
                          \wj (\zj B-A)^{-1} B \Yblock$
      \State $\Yblock := \U_{0}$ ,
    \end{algorithmic}%
  }
\end{center}
involving multiple passes of a set of initial vectors through
\eqref{eqn:contourexact} before the construction of the full subspace
using multiple moments.

Here we consider the use of full
iterations with SSM, using a linear combination of the approximate
eigenvectors as the initial vectors in the next iteration, \cite{doi:10.1260/1748-3018.7.3.249}
$\Yblock := \Xapprox R$ with
$R \in \mathbb{C}^{\mathrm{width}(\Xapprox) \times \nrhs}$.
We will call this technique the ``outer iteration.''
When the approximation of the contour integral dominates computational
expense, the inner and outer iteration for a similar size of the
subspace $\U_{0}$ will have similar cost.

Previous work \cite{Imakura2016} has considered the theoretical
filtering and convergence behaviour of the inner iteration.
This is substantially more nebulous for the outer iteration, and this work will focus on a heuristic analysis of best practices for convergence and adaptivity.  
As adaptivity relies on iterative behaviour, we focus on situations where
inner and outer iterations are truly relevant; using a smaller subspace,
and thus fewer RHS (single linear solves) per iteration.

Under a constrained subspace size, however, moment-based schemes have
been observed to stagnate over iterations before all residuals of the
computed eigenpairs have reached the prescribed tolerance, in particular
for very small tolerances.
A more flexible method might therefore allow a switch from a
multi-moment to single-moment scheme if stagnation is detected.
The single-moment scheme may be more costly, as it requires more RHS
in each linear system for the same size of subspace $\U$, but may
converge to a smaller tolerance due to the increased numerical
stability of the projected subspace.
Since this switch is typically only required when residuals are already
rather small, we do not expect to require many of the more expensive
single-moment iterations.
Furthermore, adaptive techniques such as the locking of converged
eigenpairs may reduce the cost of a more expensive single-moment
iteration by reducing the number of RHS considered.
The heuristics and advantages of these two techniques will provide the
main contribution in this work.

\section{Inner vs.\ outer iterations and multi- vs.\ single-moment}

In this section we compare the convergence behaviour of iteration
schemes using either only inner iterations or only outer iterations,
or switching from inner to outer after iteration~$8$, as well as
consider the effect of switching from multiple moments (SSM-type) to
single-moment (FEAST-type).
This motivates the adaptive strategy described in the following
section.

We considered the toy problem 
\begin{equation}
  \label{eq:ToyProblem}
  A \mathbf{x}_i = \lambda_i \mathbf{x}_i
  , \quad
  \lambda_i \in \Ilambda = \left[-1,1 \right]
	  , \qquad \mbox{where }
  A = \mathrm{diag}\left(-2.99,-2.89,...,6.91\right)
        \in \mathbb{R}^{100 \times 100}
  ,
\end{equation}
containing $\ntrue = 20$ eigenvalues in the search interval.

The unit circle was chosen as the contour of integration, and we used
Gauss--Legendre quadrature with $16$ quadrature nodes on the contour. Note that only the 8 linear systems along the upper half of the contour must be solved due to symmetry; this is the value originally defined in FEAST \cite{polizzi2009density}.
The subspace was constrained to yield a truly iterative scheme.
More precisely, the number of initial vectors after each iteration was
chosen as to give an overall subspace size
$\mo = 32 \approx 1.5 \times \ntrue$.
Testing was done in Matlab, with the initial $\Yblock$ chosen randomly using the ``twister'' generator from a standard normal
distribution, and linear systems were solved using Matlab's backslash
function.

\begin{figure}
  \begin{center}
  \includegraphics[width=\textwidth]{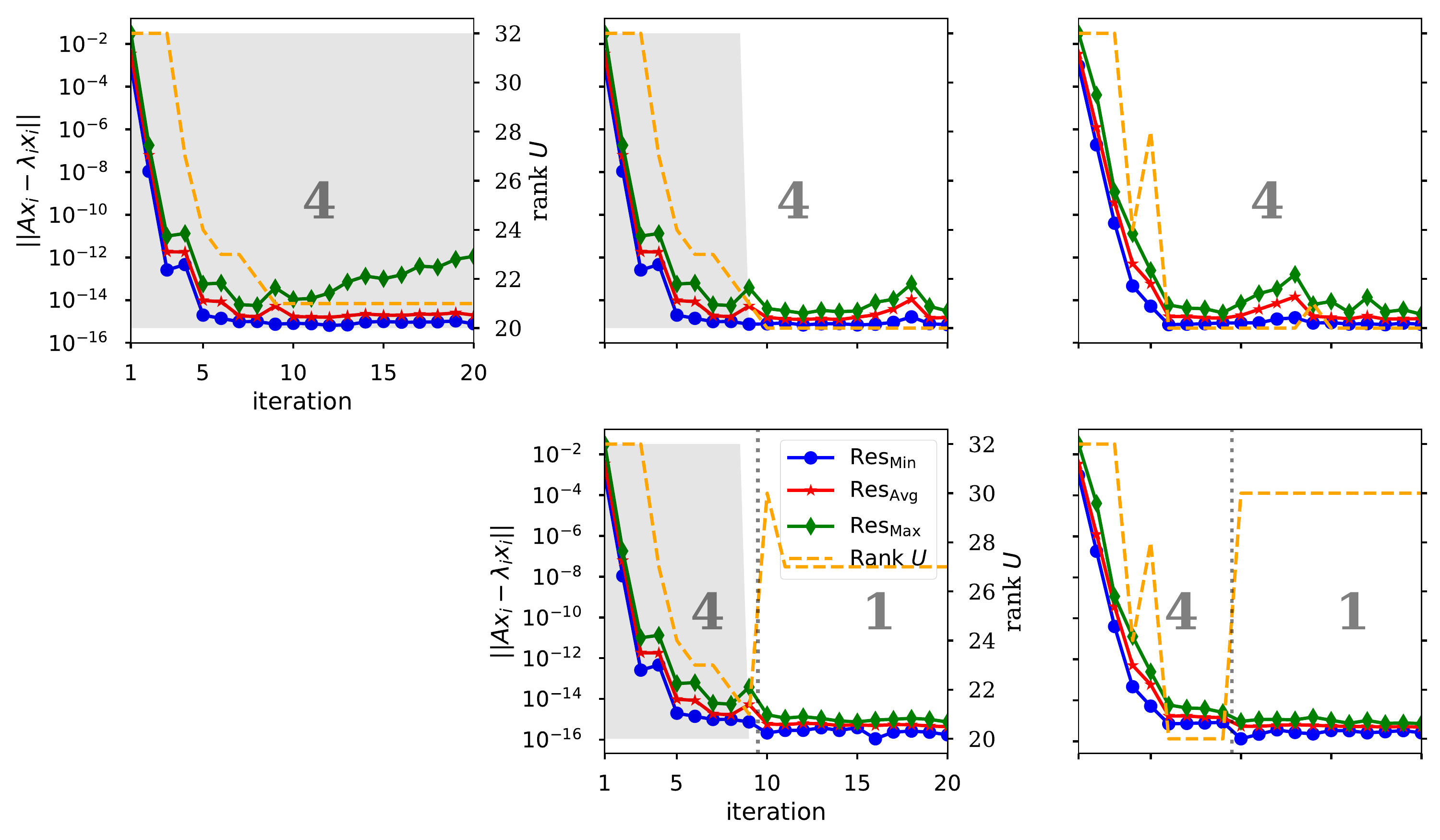}
  \end{center}
  \caption{Rank and convergence (maximum, minimum, and average of the
    $20$ smallest residuals) of the contour integration-based
    iterative schemes
    for problem \eqref{eq:ToyProblem} with initial subspace size $32$
    and $4$ moments for the first iteration. The number of moments is denoted in the background, with a 			vertical line indicating a change in number of moments.  Background hashed for inner 	iterations.
    Bottom row: Switch from SSM-type multi-moment
      (\beastm implementation, cf.\ Section~\ref{sec:BEAST})
	to FEAST-type single-moment
	(\beastc, cf.\ Section~\ref{sec:BEAST})
      after the $9$th iteration. (I.e., single moment used in 10th and subsequent iterations.)
    Top row: No such switching, multi-moment throughout.
    Left column: All inner iterations.
    Middle column: Inner iterations for the first $8$ iterations, then outer
      iterations. (I.e., initial subspace in 10th iteration is constructed from approximate eigenvectors of previous iteration.)
    Right column: All outer iterations.
    Bottom left: Left blank as inner iterations are not consistent with FEAST-type single 	moment.
  }
  \label{fig:convergenceinnerouterswitch4moments}
\end{figure}

\begin{figure}
  \begin{center}
  \includegraphics[width=\textwidth]{%
    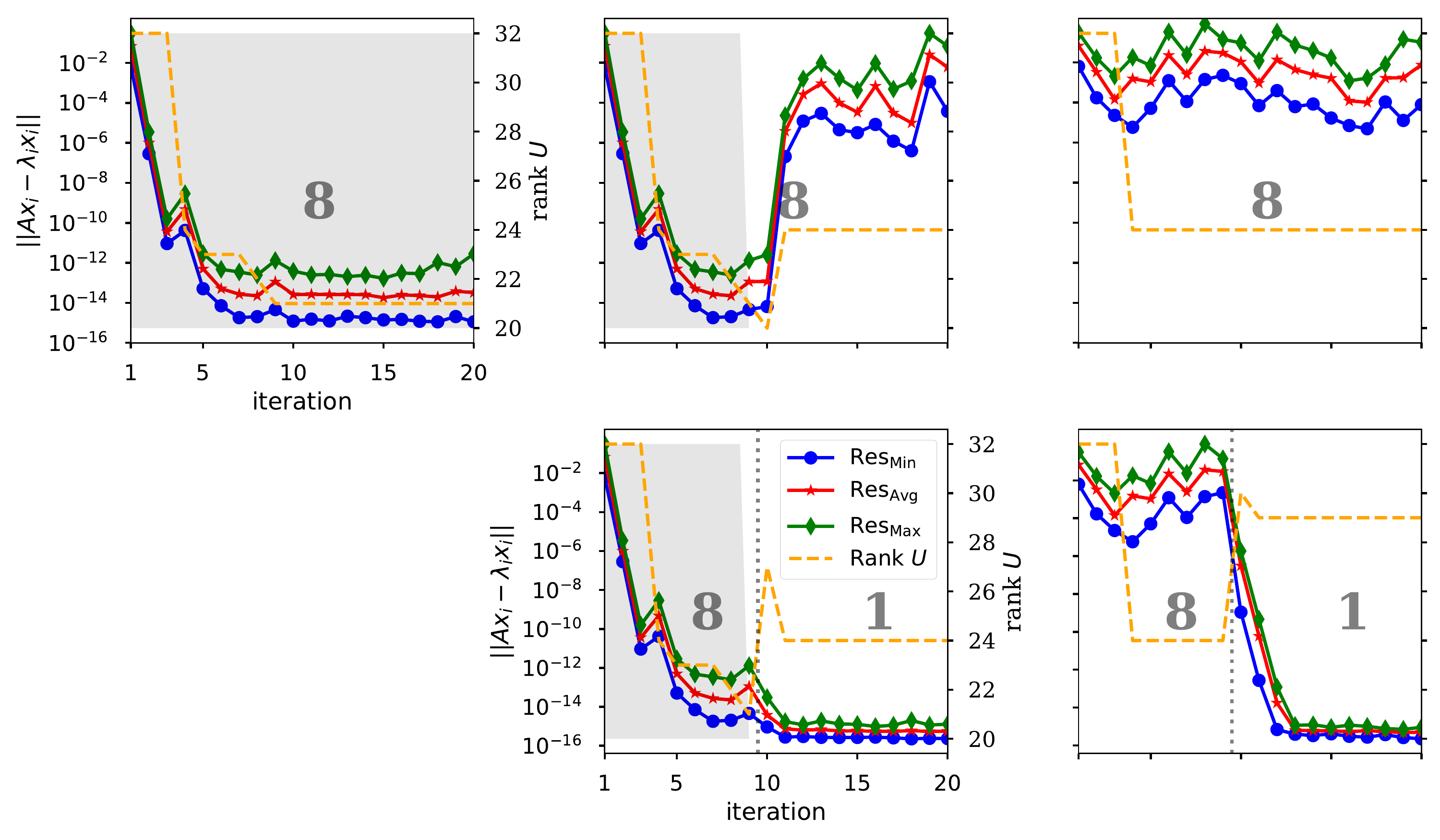}
  \caption{As Figure~\ref{fig:convergenceinnerouterswitch4moments},
    but with $8$ moments.
  }
  \label{fig:convergenceinnerouterswitch8moments}
  \end{center}
\end{figure}

The minimum, average, and maximum of the smallest $\ntrue = 20$ residuals for
the approximate eigenpairs are plotted in
Figure~\ref{fig:convergenceinnerouterswitch4moments}
(using $4$ moments and therefore $8$ columns in $\U_{0}$) and
Figure~\ref{fig:convergenceinnerouterswitch8moments}
($8$ moments, $4$ columns).
We observe different behaviour for the inner and outer iteration types
with moments, which further varies with the number of moments (and size
of $\U_{0}$).

\subsection{Inner and outer iterations}

In the top row of each figure, we observe the different behaviour for
outer and inner iterations. 
Though computation of the approximate eigenpairs is not strictly necessary for the inner iteration, we include it in this work in order to observe intermediate convergence behaviour and utilize further advantageous strategies dependant on this behaviour. 
For large problems, the cost of computing $\U$ vastly outweighs the cost of the Rayleigh-Ritz step, so performance should not be significantly affected by the additional work, which may also be offset by savings through adaptive strategies.

We see that for early iteration counts, convergence is noticeably
improved for the inner iteration compared to the outer iteration.
These effects are amplified as the starting number of moments increases, as seen in the difference between Figures \ref{fig:convergenceinnerouterswitch4moments} and \ref{fig:convergenceinnerouterswitch8moments}.
In later iterations, however, some convergence may be lost as error
creeps back in to the largest residual.

\subsection{Switching from multi- to single-moment}

In the bottom row we observe the effect of switching from a multi- to
single-moment scheme.
After the $9$th iteration, we switch from a SS--RR- to a FEAST-type
iteration, with only one moment.
As the overall subspace size remains the same $1.5\times$ the
approximate number of eigenpairs, this means that more right hand sides
are involved in the solution of the linear systems.
We see the increase in stability of convergence with this result, for
previous iterations of both the outer and inner type.

\section{The \beast framework}%
  \label{sec:BEAST}

\beast provides a framework for subspace iteration with Rayleigh--Ritz
extraction of eigenvalues and vectors
\cite{galgon2014improving,galgon2015parallel}.
Testing and development are often done with a Matlab version of \beast,
and a C version with multiple levels of parallelization is available
for performance-critical runs.

\beast provides two contour-based ways to construct the subspace $\U$;
our SS--RR- and FEAST-type implementations are referred to as \beastm
and \beastc, respectively.
(For standard eigenvalue problems
$A\xexact = \xexact \ewexact$ we also provide \beastp, which uses a
polynomial filter: $\U = p( A ) \Yblock$.
This component is not considered in the present paper.)
The framework contains algorithmic strategies discussed in previous
works for improving the foundational methods, such as the locking of
converged eigenvectors \cite{NLA:NLA2124}, as well as the components
integral to the reference algorithms, such as the orthogonalization of
the subspace for \beastm.

\begin{figure}
	\includegraphics[width=\textwidth,height=7cm]{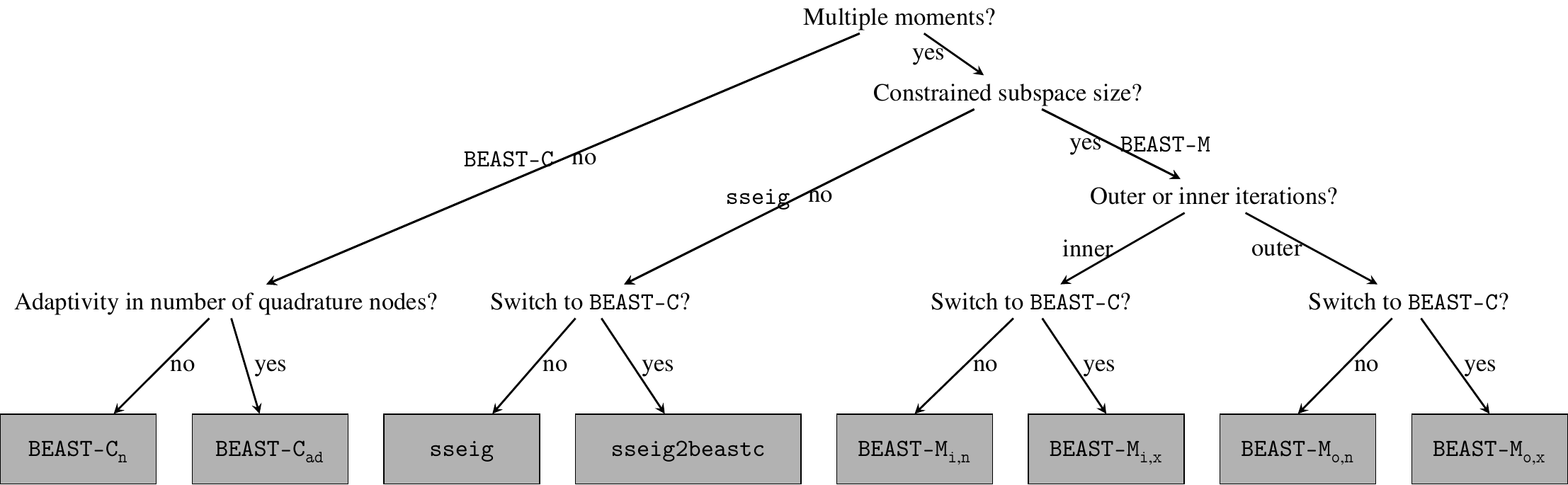}
	  \caption{Overview of all computation schemes and choices considered in section \ref{sec:numericalexperiments}.
	}
    \label{fig:methods}
\end{figure}

One particular feature of \beast is that it allows the method (or its
parameters) to be changed from one iteration to the next, and some
of the resulting overall strategies are discussed in the remainder of this section.
We first explain the naming conventions.
\begin{itemize}
  \item \beastm and \beastc refer to the method used in a given
    iteration.
  \item For \beastm, we can do inner or outer iterations
    (denoted as \beastmin and \beastmout, resp.),
    and we may or may not allow switching to \beastc whenever
    stagnation of convergence is detected (\beastmx and \beastmn,
    resp.).
    Thus, $\beastmxout$ starts with multi-moment outer iterations and
    may switch to single-moment upon stagnation.

    Stagnation was measured in the ``drop rate'' ratio of the smallest
    residual among the not-yet-converged eigenpairs,
    $r_{\snc}
     = \max_{j} \{ \| A\xapprox_{j}-B\xapprox_{j}\ewapprox_{j} \| :
                   \mbox{this residual is above tolerance} \, \}$,
    since the previous iteration; that is, $r^i_{\snc}/r^{i-1}_{\snc}$.
    Whenever this ratio rose above the threshold of $0.01$, stagnation
    was determined to have set in.  
  \item For \beastc iterations, the variant \beastcad adjusts the number
    of quadrature nodes $q$ (and the nodes and
    coefficients) over the iterations, whereas \beastcn does not.
\end{itemize}
Figure~\ref{fig:methods} shows the different choices for the methods, as well as SSM-based methods included in section \ref{sec:numericalexperiments}, 
and Algorithm~\ref{alg:beastfull} gives an overview of the computations
done in each \beast case.

\begin{algorithm}
  \caption{\beast framework, including \beastm and \beastc algorithmic
    variants.
    The option of adaptive quadrature node choice for \beastc is shown
    as \beastcad.
    The option of adaptive method switching for \beastm is shown as
    \beastmx.
    The option of inner or outer iterations for \beastm are shown as
    \beastmout and \beastmin.
  }
  \label{alg:beastfull}
  \begin{algorithmic}[1]
    \Input Matrix pair $A,B$ of size $\Nmat \times \Nmat$ 
    \Inputcont $\nexpect$ estimate for number of eigenpairs in interval
    \Inputcont $\quadpts$ quadrature nodes and coefficients
      ($\z_j, \omega_j$)
    \Inputcont $\Yblock \in \mathbb{C}^{\Nmat \times \nrhs}$
    \While{not converged}
      \If{in \beastcad mode}
        \State Choose $\quadpts, \z_j, \omega_j$
          \Comment{otherwise use parameters as provided with inputs}
      \EndIf
      \State Construct subspace $\U$ according to
        \eqref{eqn:momentssubspace} (\beastm mode) or
        \eqref{eqn:contourexact} (\beastc mode)
      \If{in \beastmin mode}
        \State Set $\Yblock:=U_0$ \Comment{inner iteration}
	  \EndIf      
        \State Compute singular values of $\U$ 
        \If{in \beastm mode}
          \State Orthogonalize $\U$
        \EndIf
        \State Resize subspace
        \State Rayleigh--Ritz extraction of eigenvalues $\Lambda$ and
          eigenvectors $\Xapprox$ according to \eqref{eq:RR}
        \State Orthogonalize $\Xapprox$ against locked eigenpairs,
          lock newly converged eigenpairs
        \If{in \beastmx mode and stagnation has been detected}
        	\State switch to \beastc
      	\EndIf
	  \If{in \beastmout mode}
          \State Set $\Yblock:=\Xapprox R$ with
            $R \in \mathbb{C}^{\mo \times \nrhs}$
	  \ElsIf{in \beastc mode}
	    \State Set $\Yblock := \Xapprox$
	  \EndIf
    \EndWhile
  \end{algorithmic}
\end{algorithm}

\section{Numerical experiments}
\label{sec:numericalexperiments}

Numerical results are shown for $37$ test problems as detailed in Table \ref{tab:problems}. The problems arise from graphene modelling \cite{RevModPhys.81.109}, and the SuiteSparse Matrix Collection \cite{uflorida}.

All contour-based schemes require selecting a quadrature rule.
Gauss--Legendre is a typical choice for FEAST, and the trapezoidal rule is
a typical choice for SSM.
The reasoning for these standards becomes clear when comparing the
$\RHSovl$ counts over several test problems.
We used an elliptical contour with eccentricity $0.1$ and 16 quadrature nodes.
Testing \beastmxout and \beastmnin with both the large subspace size
typical for SSM (dimension of overall subspace four times larger than
the expected number of eigenvalues) and the constrained size typical for
FEAST (factor $1.5$), in both cases, the average number of $\RHSovl$ was
similar or better for Gauss--Legendre with subspace factor $1.5$, and
trapezoidal was superior with subspace factor $4$; cf.\
Table~\ref{tab:quadrule}.
The difference was especially marked in the case of the outer iteration.
Investigating the causes of this difference is left for future work.
In the following tests, we used Gauss--Legendre for \beast and the
trapezoidal rule for \sseig (a Matlab implementation of SS--RR \cite{doi:10.1260/1748-3018.7.3.249}),
in both cases for an elliptical contour with eccentricity $0.1$ and 16 quadrature nodes on the whole contour. 

\begin{table}
	\centering
   	\caption{Sizes, number of
     nonzeros per row (nnzr; rounded), smallest ($\gmin$) and largest
     ($\gmax$) eigenvalues (rounded to three significant digits), and
     the two search intervals (with eigenvalue counts $\ntrue$) for 9 test
     matrices from graphene modelling and 10 test matrices from the
     GHS\_indef (laser, linverse, brainpc2), PARSEC (SiH4, Si5H12, SiO), 			 ACUSIM (Pres\_Poisson), Boeing (bcsstk37), DIMACS10
     (rgg\_n\_2\_15\_s0), and Andrews (Andrews) groups of the
     SuiteSparse Matrix Collection. \cite{uflorida}
     For each problem, the test matrix and
     interval were taken as $A$ and $\Ilambda$ in \eqref{eqn:genevp},
	 $B=I$. The problems are also problem numbers 1:25, 27:34, 36:39 in 			 Table~1 of Galgon et al.
	\cite{galgon2014improving}.}
	\label{tab:problems}
	\begin {tabular}{llllllllll}%
	\rowcolor {lightgray}Name&Size&nnzr&$\left [ \gmin \textrm {,} \gmax 			\right ]$&Interval A&$\ntrue $&No.&Interval B&$\ntrue $&No.\\%
	laser&3002&3&$\left [-1.10,4.25\right ]$&$\left [-0.100,0.357\right ]			$&307&1&$\left [1.0000,4.2389\right ]$&304&NA\\%
	\rowcolor [gray]{0.9}SiH4&5041&34&$\left [-0.996,36.8\right ]$&$\left 			[25.0,28.4\right ]$&301&2&$\left [15.3,16.4\right ]$&290&3\\%
	linverse&11999&8&$\left [-4.70,15.5\right ]$&$\left [0.00,0.62\right ]			$&308&4&$\left [2.62,2.77\right ]$&304&5\\%
	\rowcolor [gray]{0.9}Pres\_Poisson&14822&48&$\left [1.28 e{\text -} 			5,26.0\right ]$&$\left [3.7,10.0\right ]$&302&6&$\left [1.000,1.182\right 	]$&300&7\\%
	Si5H12&19896&37&$\left [-0.996,58.6\right ]$&$\left [24.2,24.7\right ]			$&320&8&$\left [41,42\right ]$&274&9\\%
	\rowcolor [gray]{0.9}bcsstk37&25503&45&$\left [-7.04 e {\text -} 				5,8.41e7\right ]$&$\left [8.0e5,9.3e5\right ]$&297&10&$\left 					[1.15e7,1.30e7\right ]$&305&11\\%
	brainpc2&27607&6&$\left [-2000,4460\right ]$&$\left [275,345\right ]			$&313&12&$\left [1900,1920\right ]$&309&13\\%
	\rowcolor [gray]{0.9}rgg\_n\_2\_15\_s0&32768&10&$\left [-5.12,17.4\right 		]$&$\left [-2.00,-1.95\right ]$&282&14&$\left [5.0,5.5\right ]$&337&15\\%
	SiO&33401&39&$\left [-1.67,84.3\right ]$&$\left [32.0,32.4\right ]				$&316&16&$\left [57.0,57.8\right ]$&283&17\\%
	\rowcolor [gray]{0.9}Andrews&60000&13&$\left [3.64e {\text -} 					16,36.5\right ]$&$\left [21.0,21.4\right ]$&300&18&$\left 						[11.20,11.26\right ]$&298&19\\%
	GraI-1k&1152&13&$\left [-3.43,2.78\right ]$&$\left [-0.7575,1.1025\right 		]$&301&20&$\left [0.42,1.58\right ]$&301&21\\%
	\rowcolor [gray]{0.9}GraI-11k&11604&13&$\left [-3.43,2.78\right ]$&$\left 	[-0.0475,0.3925\right ]$&298&22&$\left [0.935,1.065\right ]$&289&23\\%
	GraI-119k&119908&13&$\left [-3.43,2.78\right ]$&$\left 							[0.1375,0.2075\right ]$&306&24&$\left [0.9957,1.0043\right ]$&304&25\\%
	\rowcolor [gray]{0.9}GraII-1k&1152&13&$\left [-3.43,2.78\right ]$&$\left 		[-0.7575,1.1025\right ]$&292&26&$\left [0.42,1.58\right ]$&304&27\\%
	GraII-11k&11604&13&$\left [-3.43,2.79\right ]$&$\left 							[-0.1375,0.4825\right ]$&299&28&$\left [0.949,1.051\right ]$&299&29\\%
	\rowcolor [gray]{0.9}GraII-119k&119908&13&$\left [-3.43,2.79\right ]$&$			\left [0.0975,0.2475\right ]$&313&30&$\left [0.9956,1.0044\right ]				$&303&31\\%
	GraIII-1k&1152&12&$\left [-3.35,2.73\right ]$&$\left 							[-0.7575,1.1025\right ]$&300&32&$\left [0.42,1.58\right ]$&331&33\\%
	\rowcolor [gray]{0.9}GraIII-11k&11604&12&$\left [-3.35,2.73\right ]$&$			\left [-0.0975,0.4425\right ]$&305&34&$\left [0.952,1.048\right ]				$&319&35\\%
	GraIII-119k&119908&13&$\left [-3.43,2.78\right ]$&$\left 						[0.1395,0.2055\right ]$&310&36&$\left [0.996,1.004\right ]$&311&37\\%
\end{tabular}%
\end{table}

\begin{table}
	\centering
	\begin{minipage}{0.6\linewidth}
  	\caption{Average $\RHSovl$
	with different quadrature rules for 31 of the 37 problems listed in Table 	\ref{tab:problems}.  Problems for which all not schemes found all 				eigenpairs (7,14,20,26,34) are not included in the average.}
  	\label{tab:quadrule}
	\begin{tabular}{|l|l|r|r|}
		\toprule
		Iteration Type & Quadrature Rule &  Subspace Factor &  Average $				\RHSovl$ \\
		\midrule
		   \beastmxin &  Gauss Legendre &              1.5 &         3587 \\
		   \beastmxin &     Trapezoidal &              1.5 &         3514 \\    
		   \beastmxin &  Gauss Legendre &              4.0 &         5229 \\
		   \beastmxin &     Trapezoidal &              4.0 &         4545 \\
		   \beastmxout &  Gauss Legendre &              1.5 &         3201 \\
		   \beastmxout &     Trapezoidal &              1.5 &         3723 \\
		   \beastmxout &  Gauss Legendre &              4.0 &         5218 \\
		   \beastmxout &     Trapezoidal &              4.0 &         4366 \\
		\bottomrule
	\end{tabular}
\end{minipage}%
\qquad
\begin{minipage}{0.35\linewidth}
	\caption{Number of problems for each method for which not all
    eigenpairs were found (from a total of $37$ problems).}
  	\label{tab:missedcount}
	\begin{tabular}{|l|r|}
		\toprule
		         Solver &  Total failed \\
		\midrule
		      \beastcad &             7 \\
		         \sseig &             7 \\
		 \sseigtobeastc &             1 \\
		     \beastmnin &             8 \\
		     \beastmxin &             2 \\
		    \beastmnout &             8 \\
		    \beastmxout &             1 \\
		\bottomrule
	\end{tabular}
\end{minipage}
\end{table}

We first compared the total number of right hand sides used in all
linear system solves, $\RHSovl$, as well as all block linear solves,
$\BLSovl$, over the iterations required to reach a tolerance of
$10^{-13}$.
Testing was done in Matlab, with the initial $\Yblock$
chosen randomly and linear systems solved using Matlab's backslash
function.
In \beast, orthogonalization and the computation of singular values was
done in one step using Matlab's \texttt{svd} function.
Small eigenproblems were solved directly using the \texttt{eig}
function. 
In accordance with FEAST \cite{polizzi2009density} we chose a subspace
size $\mo \approx 1.5 \times \nexpect$, where $\nexpect$ is the
expected number of eigenpairs in the interval.
Given an estimated number of eigenpairs $\nexpect = 300$ for all test problems,
\beastm and \beastc began with overall subspace sizes $\mo$ of $452$ and
$450$, respectively.
\beastm began in all cases with $4$ moments.

We compared with the \beastc solver with an adaptive number of
quadrature nodes $q$, \beastcad, and with the \emph{best} fixed-order
method, \beastcn (i.e., the fixed $q$ leading to the minimum $\RHSovl$
among all $q$ considered, such that that all eigenpairs were found).
In the adaptive case, as described by Galgon et al. \cite{galgon2014improving}, $q$ is increased by a factor of $\sqrt{1.5}$ or at least by 1 if the smallest non-converged residual decreases by a factor of $\left( 10^{-1.5}, 10^{-0.75} \right]$ and by a factor of $1.5$, and by at least 1, if the decrease factor is greater than $10^{-0.75}$.
A comparison was also made to \sseig, which was also run with $4$
moments, and the subspace was allowed to grow expansively, starting
with $\mathrm{width}( U_{0} ) = 16$ right hand sides and growing to a
maximum of $256$ for the problems considered, for an overall subspace
size of $\mo=1024$. 
The final comparison was with \sseigtobeastc, switching from \sseig
after a single outer iteration to \beastc, using all approximate
eigenvectors as the initial vectors $Y$ for \beastc.
This allowed us to test the scheme-switching heuristic in the context
of large subspace sizes.
With a large subspace size, it is expected that \sseig and
\sseigtobeastc will converge quickly, i.e., with just a few iterations.

\begin{figure}
    \centering
    \begin{minipage}{.5\textwidth}
        \centering
  \includegraphics[width=9cm]{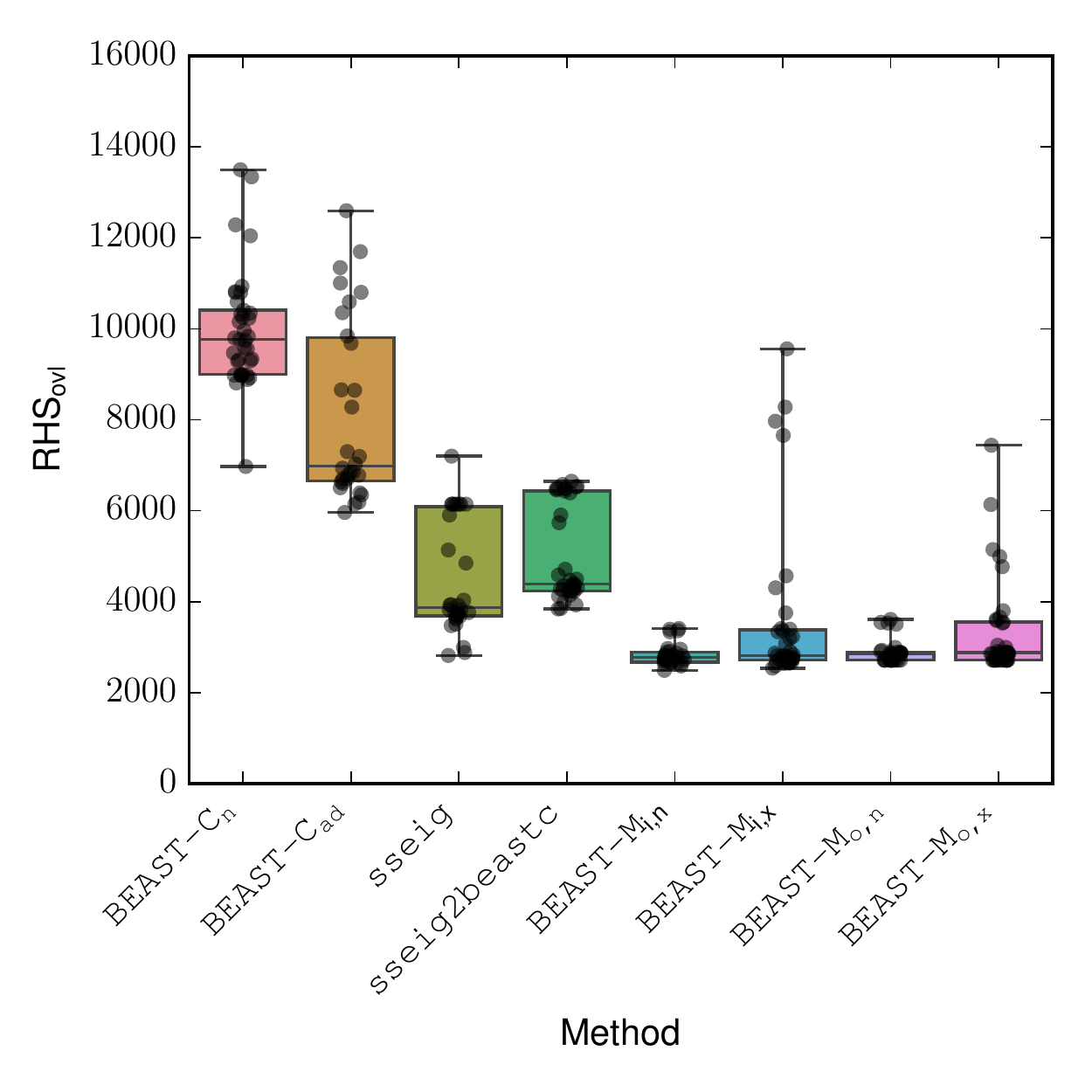}
  \label{fig:matlabSLS}
    \end{minipage}%
    \begin{minipage}{0.5\textwidth}
      \centering
  \includegraphics[width=9cm]{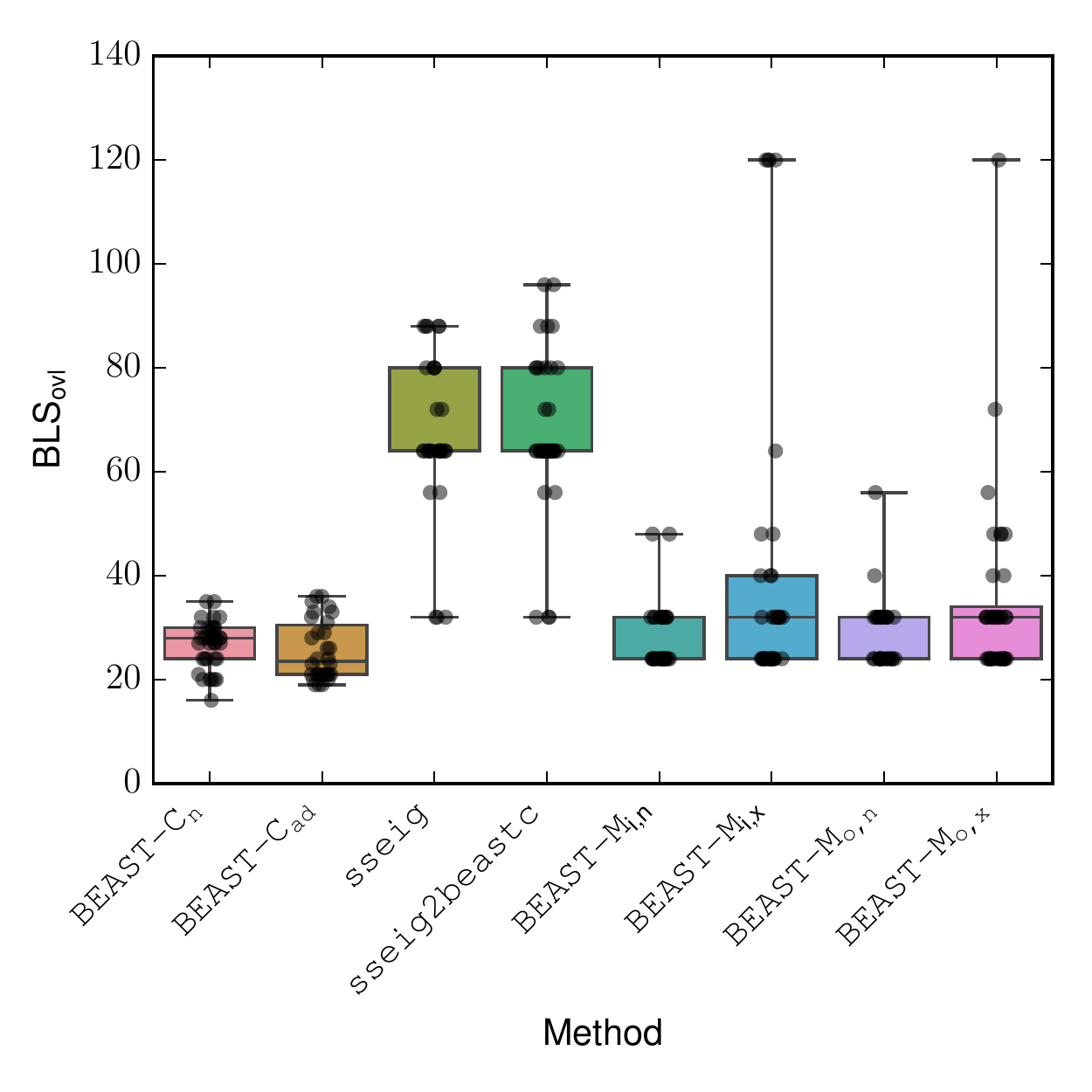}
  \label{fig:matlabBLS}
    \end{minipage}
      \caption{$\RHSovl$ (left) and $\BLSovl$ (right) for different solution 		  techniques. Values are only shown for systems where all eigenpairs were 	  found with a given method. Each box extends from the first to third 			  quartiles of the respective data, with the median (second quartile) 			  shown as a horizontal black line in the box. The whiskers cover the 			  entire range of observations.
      Distinct counts for each problem are visible as translucent circles. 			  Full numerics are given in the appendix.
  }
  \label{fig:matlabLS}
\end{figure}

We measured success in terms of the number of linear solves (single
or block) required for all eigenpairs to converge (efficiency), as well
as the number of times the method found all eigenpairs inside the
interval (robustness).
The respective results are summarized in Figure \ref{fig:matlabLS}, as well as Table \ref{tab:missedcount}.
We point out that the methods were deemed to have ``failed'' when not
all eigenvalues inside the interval were found to the desired tolerance.
This may be an overly strict measure, especially for \sseig, which we
observed often had a few eigenvalues slightly above the residual
tolerance at completion.  

We observe that all methods obtaining the minimum, or near minimum
number of single linear solves, $\RHSovl$, included multiple moments,
and the difference in $\RHSovl$ between methods involving and not
involving moments is substantial.
We also observe that methods involving a switch to \beastc behaved more
robustly, with fewer failures to find the minimum number of linear
solves, $\RHSovl$.
We see that both \beastmx and \beastmn outperformed \sseig and \beastc,
with \beastmn and \beastmx having the same minimum $\RHSovl$ in most
cases.

When switching was allowed, \beastmxout and \beastmxin behaved robustly,
failing to find all eigenvalues for only one and two problems respectively.
In comparison, \beastcad , \beastmn, and \sseig featured an observably
higher number of cases where not all eigenpairs were found to the desired
tolerance.
For \beast, even if this were an acceptable behaviour, a failed convergence check (performed using singular values of $\U$) will result in continued iterations until all eigenpairs have reached the desired tolerance or some maximum number of iterations is reached. For stagnant eigenpairs this may be a considerable expense.
We also observe the increased $\BLSovl$ for \beastm compared to \beastc
and \sseig, due to a higher number of iterations.
We predict, and demonstrate, in the following tests, that the reduction
of $\RHSovl$ is enough to offset this cost, but acknowledge that this
depends on the linear solver used, the computational set-up, including
the distribution of work and memory, and the problem itself. 

\begin{figure}
  \centering
  \begin{minipage}{.5\textwidth}
    \centering
    \includegraphics[width=\linewidth]{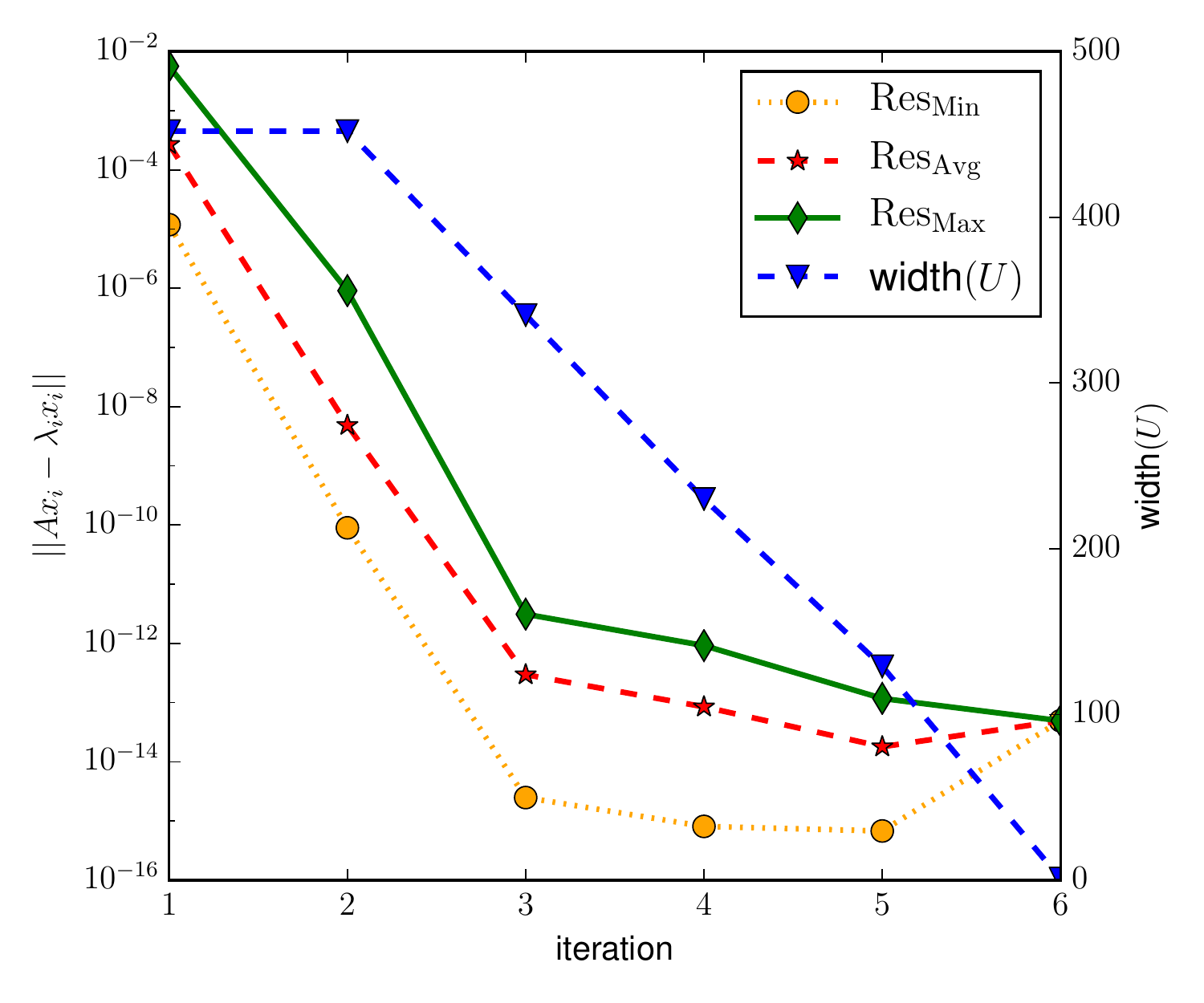}
    \subcaption{\beastmxin}
    \label{fig:beastmxin}
  \end{minipage}%
  \begin{minipage}{0.5\textwidth}
    \centering
    \includegraphics[width=\linewidth]{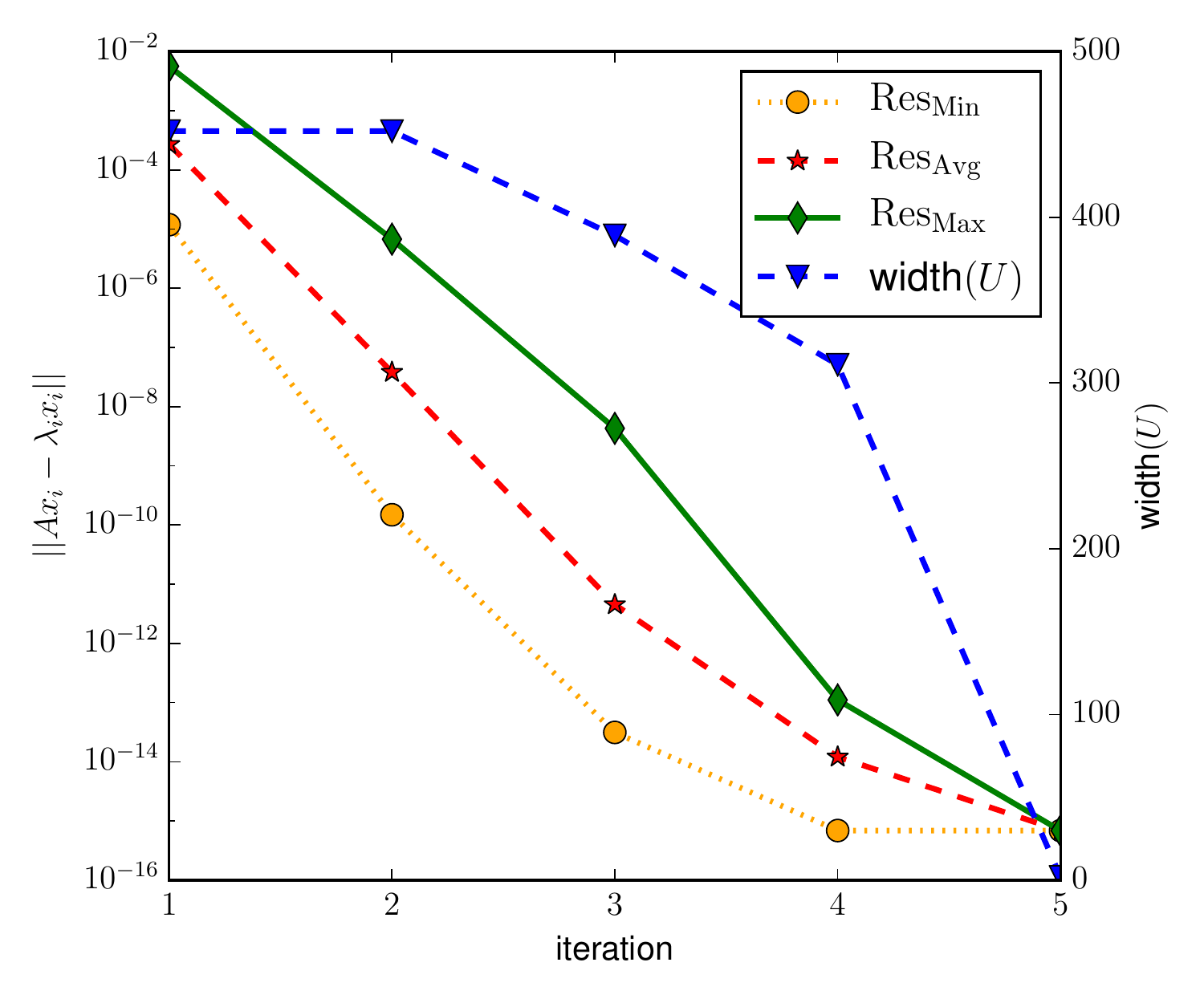}
    \subcaption{\beastmxout}
    \label{fig:beastmxout}
  \end{minipage}
  \caption{Convergence and subspace size $\mo$ with
    \beastmxin and \beastmxout for Problem 4, Table 1 in Galgon et al.
    \cite{galgon2014improving}.
    }
  \label{fig:innnerouteriteration}
\end{figure}

Taking a closer look at \beastmout and \beastmin, we see similar
behaviour between the two methods, with a few larger outliers in
$\RHSovl$ for \beastmxin.
Based on our results from the previous section, we might have expected
noticeably better performance from the inner iteration as opposed to the
outer.
In fact, the initial improvement in convergence for the inner iteration
is a significant reason for these outliers.
We consider the numerical example shown in
Figure~\ref{fig:innnerouteriteration}.
After convergence begins, eigenpairs are locked, and the subspace
shrinks; perhaps overly much, as convergence subsequently slows down,
even after a switch to \beastc.
Eventually, an additional iteration is required for \beastmxin, which
started out with much faster convergence than \beastmxout. 

Out of the various methods considered so far, we choose to highlight
\beastmxout for its combination of robustness, adaptability, and low
cost.
We now turn to numerical results showing the possible savings in cost
available with this scheme.%

We also considered how a reduction in $\RHSovl$ corresponds to savings
in time.
We compared parallelized, distributed memory implementations of \beastc
(specifically \beastcn) and \beastmxout for two larger systems.
The scalable solver \strumpack \cite{strumpack} was used to solve all
linear systems directly, and the kernel library \ghost
\cite{kreutzer2017ghost} provided fast matrix operations as a
background library.
Testing was performed with $32$ nodes of the Emmy HPC cluster at
Friedrich-Alexander-Universit{\"a}t Erlangen--N{\"u}rnberg.
The eigenproblems are detailed in Table~\ref{tab:problemstwo}.
They were solved to a tolerance of $10^{-12}$, using $16$ quadrature
nodes and Gauss--Legendre quadrature on a circular
contour.
The number of columns of $\U$ was set to $1.5 \times \nexpect$, and
locking was enabled.
\beastmxout began with $4$ moments and was allowed to switch to \beastc
when stagnation was detected.

\begin{table}[!h]
  \centering
  \caption{Sizes, average number of nonzeros per row (nnzr; rounded),  Approximate smallest ($\gmin$) and largest
    ($\gmax$) eigenvalues, search interval,
    starting approximate eigenvalue count $\nexpect$, and number of eigenvalues found $\hat{\ntrue}$
    for test matrices from graphene \cite{RevModPhys.81.109} and topological insulator modelling \cite{RevModPhys.82.3045}.
    Related matrices may be generated from the \scamac library. \cite{scamac} The spectral properties of matrices of these forms with regard to computation have been previously studied. \cite{ISI:000384077600012}
    For each problem, the test matrix and interval were taken as $A$ and
    $\Ilambda$ in \eqref{eqn:genevp}, $B=I$.
  }
  \label{tab:problemstwo}
	\begin {tabular}{llcllll}%
		\rowcolor {lightgray}Name&Size&nnzr&$\left [ \gmin \textrm {,} \gmax 			\right ]$&Interval&$\nexpect $&$\hat {\ntrue }$\\\hline %
		Graph16M&16000000&4&$\left [ -3.0,3.0 \right ]$&$\left [ 						-0.0025,0.0025 \right ]$&320&318\\\hline %
		Topi1M&1638400&12&$\left [ -4.8,4.8 \right ]$&$\left [ -0.06,0.06 				\right ]$&120&116\\\hline %
	\end {tabular}%
\end{table}

\begin{figure}
    {\includegraphics[width=13cm]{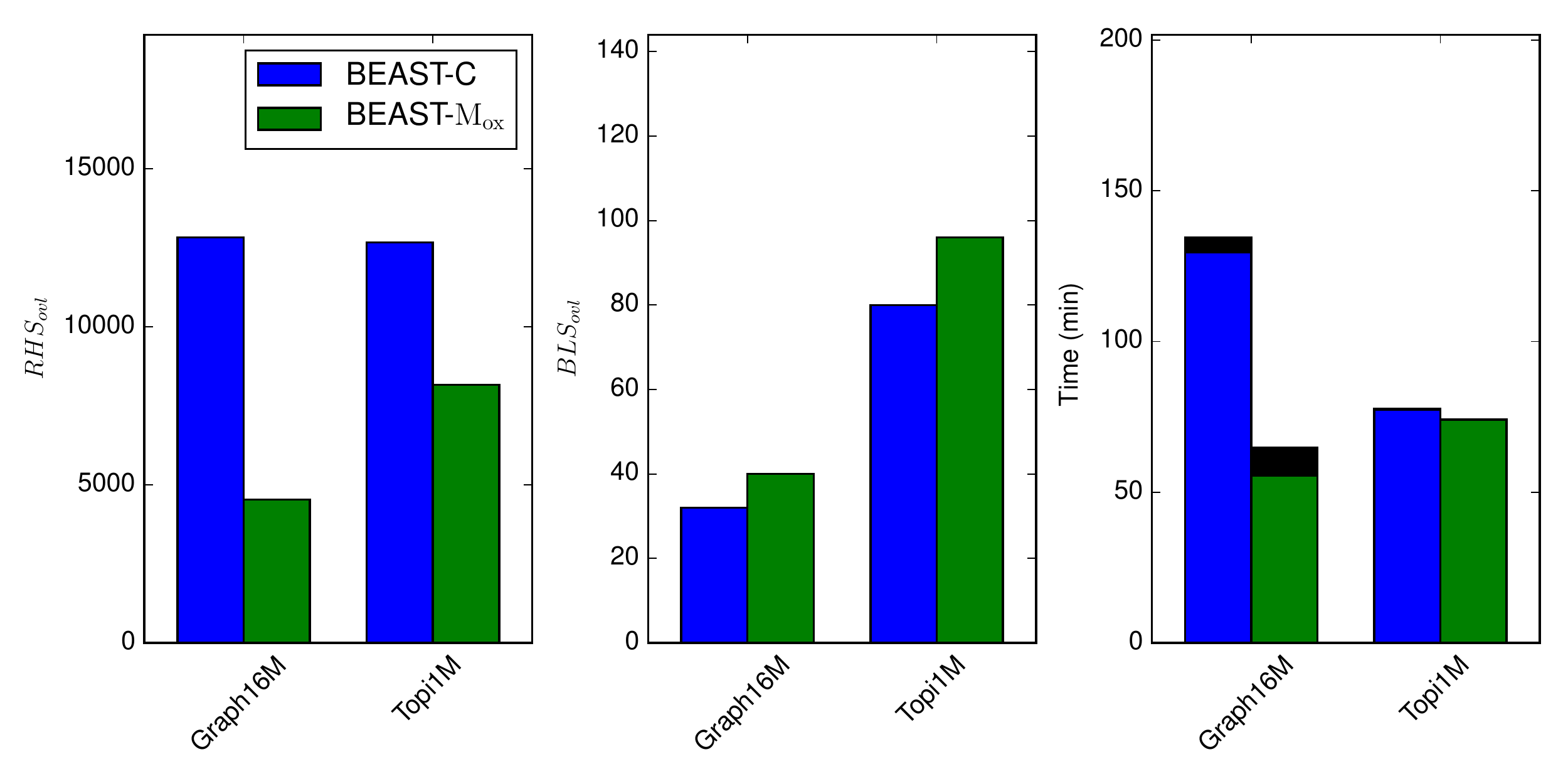}}
    \caption{Linear solve count and timings of \beastmxout vs. \beastc. In right most figure, time for building $\U$ shown in color, time for all other components shown in black.}
    \label{fig:linsolvecountandtiming}
\end{figure}

In Figure~\ref{fig:linsolvecountandtiming}, we again see the reduction in
$\RHSovl$ of \beastmxout compared to \beastc, as well as the respective times to solution.
As expected, the construction of $\U$ was by far the most time consuming
portion of \beast, due to the expense of solving linear systems.
The time reduction was more pronounced for the graphene problem than for
the topological insulator.
For the latter, the decrease in $\RHSovl$ for \beastm compared to
\beastc was not so dramatic, and at least one additional iteration was
required, as can be seen in the $\BLSovl$ count.
However, even with the additional iteration(s), \beastmxout was still
faster than \beastc.

\section{Conclusion}

Contour integration schemes are a powerful tool for interior eigenvalue
problems, but the repeated solution of large linear systems with many
right hand sides is a bottleneck in time and energy. The judicious choice of various parameters has a noticeable effect on the potential overall cost of the solver, particularily those parameters corresponding to the cost of linear system solves. Obvious parameters, such as the choice of quadrature degree, as considered in previous work \cite{NLA:NLA2124}, will have a noticeable effect on the overall cost of linear system solves per iteration, and may be chosen adaptively for an overall reduction in cost over iterations. However, this parameter, and others like it, have a direct effect on the convergence rate. We have explored the heuristical choice of the number of moments, which, while relevant to convergence, may have a less direct effect, especially for initial iterations. The choice of the number of moments is a powerful parameter for controlling the cost of linear system solves, and thus the overall eigensolver cost.
We explored further heuristics within the context of an iterative eigensolver, including best strategies for choosing the number of moments on a given iteration, and the choice of starting vectors. In a broad comparison between methods, we provide evidence that a multi-moment flexible iterative method may reduce the number of single linear solves over all iterations, and thus the potential overall cost of an eigensolver. Furthermore, the flexibility of the method ensures comparative robustness and accuracy. We provide further evidence for performance capabilities in initial larger experiments. 

Future topics of exploration include the adaptive choice of quadrature degree in combination with multiple moments, extending the work of \cite{NLA:NLA2124}. Furthermore, the capacity for performance improvement with multiple moments and reduced right hand sides clearly depends on the linear solver used. In these experiments, due to the difficulty of solving these shifted linear systems, we tested with direct solvers. In future work, the capacity for improvement dependant on the properties of the linear solver used may be investigated. Results and exploration for an iterative linear solver would be of particular interest as present problems may exceed the capacity of direct solvers.

\section*{Acknowledgements}
The present study has been supported in part by the Deutsche Forschungsgemeinschaft through the priority programme 1648 ``Software for Exascale Computing'' (SPPEXA) under the project ESSEX-II, as well as by the Japan Society for the Promotion of Science (JSPS), Grants-in-Aid for Scientific Research (Nos. 17K12690, 18H03250, 19K20280)

\bibliographystyle{IEEEtran}
\bibliography{flexiblemoments}

\appendix[Detailed results for Figure 4 \label{app1}]

\renewcommand{\arraystretch}{0.8}
\begin {longtable}{ccccccccc}%
		\caption {Matlab results $\RHSovl :\BLSovl $. Problems labelled by number as listed in Table \ref {tab:problems}. Value given as "-" when the given method did not find all eigenpairs in interval. \label {tab:matlabbeast}} \\\toprule \rowcolor {lightgray} & \multicolumn {2}{c}{\beastc }& \multicolumn {2}{c}{\texttt {sseig}} & \multicolumn {4}{c}{\beastm }\\No.&\texttt {n}&\texttt {ad}&no switch&\texttt {2beastc}&\texttt {i,n}&\texttt {i,n}&\texttt {o,n}&\texttt {o,x}\\\midrule \endhead %
	\ensuremath {1}&6972 : 16&-&5904 : 80&6392 : 80&-&8280 : 120&-&4992 : 40\\%
	\rowcolor [gray]{0.9}\ensuremath {2}&10320 : 24&-&3744 : 64&4280 : 64&2968 : 32&7968 : 120&2992 : 32&2992 : 32\\%
	\ensuremath {3}&9585 : 27&6703 : 21&7200 : 80&4120 : 64&2672 : 24&2672 : 24&2712 : 24&2712 : 24\\%
	\rowcolor [gray]{0.9}\ensuremath {4}&10234 : 28&6781 : 21&3680 : 64&4304 : 64&2776 : 24&2776 : 24&2856 : 24&2856 : 24\\%
	\ensuremath {5}&8885 : 20&6784 : 21&3920 : 64&4392 : 64&2752 : 24&2752 : 24&2712 : 24&2712 : 24\\%
	\rowcolor [gray]{0.9}\ensuremath {6}&9835 : 35&10803 : 33&6144 : 88&6576 : 96&3336 : 32&3336 : 32&3504 : 56&3584 : 48\\%
	\ensuremath {7}&9765 : 27&6697 : 21&3648 : 64&4224 : 64&2728 : 24&2728 : 24&2712 : 24&2712 : 24\\%
	\rowcolor [gray]{0.9}\ensuremath {8}&10410 : 30&-&-&4584 : 64&2808 : 24&2808 : 24&2848 : 32&2848 : 32\\%
	\ensuremath {9}&9800 : 30&6637 : 20&3472 : 64&3928 : 64&2488 : 24&2536 : 32&2888 : 24&2888 : 24\\%
	\rowcolor [gray]{0.9}\ensuremath {10}&8980 : 30&-&3712 : 64&4232 : 64&2584 : 24&2584 : 24&2712 : 24&2712 : 24\\%
	\ensuremath {11}&10818 : 27&6189 : 19&3824 : 64&4352 : 64&2808 : 32&2808 : 32&2880 : 32&2880 : 32\\%
	\rowcolor [gray]{0.9}\ensuremath {12}&10936 : 32&7194 : 24&3936 : 72&4472 : 72&2792 : 24&2792 : 24&2712 : 24&2712 : 24\\%
	\ensuremath {13}&9297 : 27&6390 : 20&3824 : 64&4384 : 64&2696 : 24&2696 : 24&2712 : 24&2712 : 24\\%
	\rowcolor [gray]{0.9}\ensuremath {14}&8980 : 30&6839 : 23&-&-&-&-&-&-\\%
	\ensuremath {15}&8995 : 20&8649 : 29&4032 : 64&4712 : 64&2904 : 32&2904 : 32&2928 : 32&3040 : 48\\%
	\rowcolor [gray]{0.9}\ensuremath {16}&10596 : 30&8657 : 29&-&4496 : 64&2832 : 24&2832 : 24&2856 : 32&2856 : 32\\%
	\ensuremath {17}&8922 : 24&6348 : 19&-&3992 : 64&2648 : 24&2648 : 24&2712 : 24&2712 : 24\\%
	\rowcolor [gray]{0.9}\ensuremath {18}&9744 : 28&7299 : 24&3920 : 72&4360 : 72&2744 : 24&2744 : 24&2712 : 24&2712 : 24\\%
	\ensuremath {19}&9963 : 27&7029 : 23&3696 : 64&4232 : 64&2728 : 48&3072 : 40&2880 : 32&2880 : 32\\%
	\rowcolor [gray]{0.9}\ensuremath {20}&10255 : 35&12598 : 36&-&6456 : 80&-&4568 : 40&-&5144 : 72\\%
	\ensuremath {21}&8810 : 20&6936 : 26&2880 : 32&3848 : 32&2664 : 24&2664 : 24&2848 : 32&2848 : 32\\%
	\rowcolor [gray]{0.9}\ensuremath {22}&9000 : 20&6848 : 21&6144 : 88&6480 : 88&2952 : 32&7656 : 120&2712 : 24&2712 : 24\\%
	\ensuremath {23}&13500 : 30&10594 : 32&6144 : 88&6520 : 96&3392 : 32&3392 : 32&-&3656 : 48\\%
	\rowcolor [gray]{0.9}\ensuremath {24}&8975 : 20&9845 : 31&3616 : 64&4256 : 64&2864 : 32&2864 : 32&2872 : 32&2872 : 32\\%
	\ensuremath {25}&9471 : 28&6505 : 21&3744 : 64&4304 : 64&-&3192 : 40&2888 : 24&2888 : 24\\%
	\rowcolor [gray]{0.9}\ensuremath {26}&12047 : 28&-&-&6432 : 80&-&4304 : 48&3528 : 40&3528 : 40\\%
	\ensuremath {27}&8976 : 24&-&2816 : 32&3840 : 32&2656 : 24&2656 : 24&-&3800 : 32\\%
	\rowcolor [gray]{0.9}\ensuremath {28}&10800 : 24&9682 : 28&6144 : 88&6488 : 88&-&3232 : 64&3544 : 32&3544 : 32\\%
	\ensuremath {29}&13344 : 32&11347 : 34&6144 : 80&6488 : 80&3360 : 32&3360 : 32&-&4768 : 56\\%
	\rowcolor [gray]{0.9}\ensuremath {30}&9338 : 21&10355 : 35&4848 : 56&5736 : 56&-&9560 : 120&-&7440 : 120\\%
	\ensuremath {31}&9562 : 28&6599 : 21&3760 : 64&4304 : 64&2672 : 24&2672 : 24&2712 : 24&2712 : 24\\%
	\rowcolor [gray]{0.9}\ensuremath {32}&12288 : 32&11698 : 33&-&6496 : 80&3408 : 32&3408 : 32&3608 : 32&3608 : 32\\%
	\ensuremath {33}&10157 : 28&6150 : 21&2992 : 32&4144 : 32&2640 : 24&2640 : 24&2904 : 32&2904 : 32\\%
	\rowcolor [gray]{0.9}\ensuremath {34}&10800 : 24&11008 : 36&6144 : 88&6536 : 88&-&-&-&6136 : 48\\%
	\ensuremath {35}&10346 : 28&8279 : 26&6144 : 80&6648 : 80&2808 : 32&2808 : 32&2840 : 32&2840 : 32\\%
	\rowcolor [gray]{0.9}\ensuremath {36}&9317 : 28&-&5136 : 64&5904 : 64&2736 : 24&2736 : 24&2712 : 24&2712 : 24\\%
	\ensuremath {37}&9289 : 28&5961 : 19&3504 : 56&4240 : 56&2888 : 48&3752 : 48&2712 : 24&2712 : 24\\\hline %
\end {longtable}%

\end{document}